\documentclass[12pt]{amsart}

\usepackage{amsmath, amssymb}
\usepackage[justification=centering]{caption}
\usepackage[frame,cmtip,arrow,matrix,line,graph,curve]{xy}
\usepackage{graphpap, color, paralist, pstricks}
\usepackage[mathscr]{eucal}
\usepackage[pdftex]{graphicx}
\usepackage[pdftex,colorlinks,backref=page,citecolor=blue]{hyperref}
\usepackage{tikz}
\usepackage{kotex}
\usepackage{array}
\usepackage{pgfplots}
\usepackage{mathrsfs}
\usetikzlibrary{arrows} 
\usepackage{tikz-cd}

\newcolumntype{P}[1]{>{\centering\arraybackslash}p{#1}}
\newcolumntype{M}[1]{>{\centering\arraybackslash}m{#1}}

\newtheorem{theorem}{Theorem}[section]
\newtheorem{proposition}[theorem]{Proposition}
\newtheorem{corollary}[theorem]{Corollary}
\newtheorem{lemma}[theorem]{Lemma}

\theoremstyle{definition}

\newtheorem{example}[theorem]{Example}

\newtheorem{remark}[theorem]{Remark}

\newcommand{\PP}{\mathbb{P}}
\newcommand{\QQ}{\mathbb{Q}}
\newcommand{\CC}{\mathbb{C}}
\newcommand{\FF}{\mathbb{F}}

\newcommand{\ZZ}{\mathbb{Z}}

\newcommand{\cO}{\mathcal{O} }
\newcommand{\cA}{\mathcal{A} }
\newcommand{\cB}{\mathcal{B} }
\newcommand{\cC}{\mathcal{C} }

\newcommand{\cG}{\mathcal{G} }

\newcommand{\cI}{\mathcal{I} }
\newcommand{\cJ}{\mathcal{J} }
\newcommand{\cK}{\mathcal{K} }
\newcommand{\cM}{\mathcal{M} }

\newcommand{\cT}{\mathcal{T} }
\newcommand{\cU}{\mathcal{U} }

\newcommand{\cZ}{\mathcal{Z} }

\renewcommand{\bar}[1]{\overline{#1}}

\newcommand{\rH}{\mathrm{H} }

\newcommand\bG{\mathbf{G}}
\newcommand\bH{\mathbf{H}}

\newcommand\bM{\mathbf{M}}

\newcommand\bS{\mathbf{S}}

\newcommand{\fsl}{\mathfrak{sl}}

\def\Sym{\mathrm{Sym} }
\def\Hom{\mathrm{Hom} }
\def\Ext{\mathrm{Ext} }

\def\Gr{\mathrm{Gr} }

\def\SL{\mathrm{SL}}
\def\GL{\mathrm{GL}}
\def\git{/\!/ }

\def\lr{\rightarrow}

\newcommand{\rank}{\mathrm{rank}\, }

\newcommand{\ses}[3]{0\lr{#1}\lr{#2}\lr{#3}\lr 0}

\begin{document}

\title[Rational curves via representation]{Rational curves in a quadric threefold\\via an $\text{SL}(2,\mathbb{C})$-representation}

\author{Kiryong Chung}
\address{Department of Mathematics Education, Kyungpook National University, 80 Daehakro, Bukgu, Daegu, 41566, Republic of Korea}
\email{krchung@knu.ac.kr}

\author{Sukmoon Huh}
\address{Department of Mathematics, Sungkyunkwan University, 2066 Seobu-ro, Suwon, 16419, Republic of Korea}
\email{sukmoonh@skku.edu}

\author{Sang-Bum Yoo}
\address{Department of Mathematics Education, Gongju National University of Education, 27 Ungjin-ro, Gongju-si, Chungcheongnam-do, 32553, Republic of Korea}
\email{sbyoo@gjue.ac.kr}


\begin{abstract}
In this paper, we regard the smooth quadric threefold $Q_{3}$ as Lagrangian Grassmannian and search for fixed rational curves of low degree in $Q_{3}$ with respect to a torus action, which is the maximal subgroup of the special linear group $\text{SL}(2,\mathbb{C})$. Most of them are confirmations of very well-known facts. If the degree of a rational curve is $3$, it is confirmed using the Lagrangian's geometric properties that the moduli space of twisted cubic curves in $Q_3$ has a specific projective bundle structure. From this, we can immediately obtain the cohomology ring of the moduli space.
\end{abstract}
\maketitle


\section{Introduction}

\subsection{Motivation and results}

By definition, a \emph{Fano threefold} $X$ is a smooth projective variety whose anti-canonical divisor $-K_X$ is ample. Such threefolds $X$ with Picard number one and second Betti number $b_{2}(X)=1$ are classified to be one of the following:
\[
\PP^3, \quad Q_3, \quad V_5, \quad \text{and}\quad  V_{22}.
\]
It is well known that the first three $\PP^3$, $Q_3$ and $V_5$ are rigid but the last one $V_{22}$ forms a $6$-dimensional family; see the excellent paper \cite[Section 5.3]{Don08}. All of these varieties have the same cohomology groups as $\PP^3$ does. Furthermore, the $\SL(2, \CC)$-orbit descriptions of these Fano threefolds are well studied in the papers \cite{MU83, AF93}.

The reason for writing this paper is to fill in the missing list for the space of rational curves. The moduli space of rational curves of degree $3$ or less in the Fano variety has been studied very diversely (\cite{MU83, PS85, EPS87, Ili94, Sch01, VX02, KS04, CK11, CHK12, San14}). In the case of $\PP^3$, starting from the analysis of the Hilbert scheme \cite{PS85}, it served as a key example of enumerative geometry in the 80s and 90s (\cite{ES96}) and created many variant spaces related with birational geometry (\cite{VX02, EPS87}). In the case of quintic del Pezzo threefold, the description of the space of rational curves was possible from the Grassmannian geometry of lines (\cite{Ili94, San14, CHL18}). Lastly, in the case of Mukai variety, its construction is very fantastic (\cite{MU83, Sch01}), and it plays a very decisive role in relation to the existence of the K{\"a}hler-Einstein metric on Fano manifold. For specific conclusions, refer to the papers \cite{Don08} and \cite{CS18}. In this paper, we examine the moduli space of twisted cubic curves in a quadric threefold. The authors think that a description of its moduli space is just not written down in the literature because it is well-known and easy to experts. However, it seems necessary to record global geometry and cohomological ring in order to serve as concrete examples required by enumerative geometry (\cite{Cao19a, CMT18, CLW21}). Now let us introduce the main contents of this paper. Let $X$ be a projective variety with a fixed embedding $X\subset \PP^r$.
\begin{itemize}
\item Let $\bS_{d}(X)$ be the moduli space of \emph{stable} sheaves $F$ on $X$ with Hilbert polynomial $\chi(F(m))=dm+1$.
\item Let $\bH_{d}(X)$ be the Hilbert scheme of curves $C$ in $X$ with $\chi(\cO_C(m))=dm+1$.
\end{itemize}
In this paper, we focus on the case that $X=Q_3$ is a quadric threefold. $Q_3$ is homogeneous and thus one can apply the main results of \cite{CK11,CHK12}. Specially, $\bS_{d}(Q_3)$ is smooth for $d\leq 3$. Also the quadric threefold $Q_3$ does not contain any plane and so $\mathbf{S}_3(Q_3)$ is isomorphic to $\mathbf{H}_3(Q_3)$. We find the torus invariant rational curves of the lower degrees and extend it into the cubic curve case. It turns out that the torus fixed loci are isolated until degree $2$, but not in degree $3$.
\begin{proposition}[\protect{Proposition \ref{mainprop1} and Proposition \ref{mainprop2}}]
The torus invariant cubic curves in $Q_3$ consist of isolated ones ($32$) and two connected components (isomorphic to $\PP^1$).
\end{proposition}
Through the computation of torus invariant cubic curves, one realizes that the cubic curves in the hyperplane section of $Q_3$ need to be ordered depending on its linear class. By using the universal family of the Fano scheme of lines in $Q_3$, we obtain a global description of $\bS_3(Q_3)$ as follows.
\begin{theorem}[\protect{Proposition \ref{mainprop4}}]\label{mainthm}
The moduli space $\bS_3(Q_3)$ of stable sheaves on $Q_3$ is isomorphic to a projectivized rank $6$ bundle over $\Gr(2,4)$.
\end{theorem}
The key point of the proof of Theorem \ref{mainthm} is to give the algebraic meaning of the ruling line of the hyperplane section of $Q_3$. From Theorem \ref{mainthm}, one can easily obtain the cohomology ring of $\bS_3(Q_3)$ (Corollary \ref{maincor}).

The method of the proof of Theorem \ref{mainthm} is a very natural approach to study the space of curves in Fano threefolds. For example, the moduli space of rational quartic curves in $Q_3$ is birational to the space of cubic curves in $\PP^3$ by using Gauss map or tangent developable surfaces (\cite[the proof of Proposition 4.17]{FGP19}). The first named author is willing to study of birational relations among these moduli spaces in the following paper.
\subsection{Organization of the paper} In Section \ref{collect}, we collect the well-known facts for finding torus invariant curves. In Section \ref{bb}, we apply the Bia\l ynicki-Birula theorem to the moduli space $\bS_{d}(Q_3)$ for $d\leq 3$. Lastly, in Section \ref{sub:defq}, we describe the moduli space $\bS_3(Q_3)$ as a projective bundle.

\subsection*{Notations and conventions}
\begin{itemize}
\item Let us denote by $\Gr(k,n)$ the Grassmannian variety parameterizing $k$-dimensional subspaces in a fixed vector space $V$ with $\dim V=n$.
\item When no confusion can arise, we do not distinguish the moduli point $[x]\in \cM$ from the object $x$ parameterized by $[x]$.
\item The set of fixed points of $X$ is denoted by $X^{\CC^*}$ under the $\CC^*$-action.
\end{itemize}

\subsection*{Acknowledgements}
The authors gratefully acknowledge the many helpful suggestions and comments of Jeong-Seop Kim, Yeongrak Kim, Wanseok Lee, Kyeong-Dong Park, Joonyeong Won during the preparation of the paper. A part of the paper has been initiated in the workshop (held in Jinju, Korea, Feb. 2022) aimed for finding research topics through arXiv preprint survey. K. Chung was supported by the National Research Foundation of Korea (NRF-2021R1I1A3045360, NRF-2022M3C1C8094326 and NRF-2022M3H3A1098237). S. Huh was supported by the National Research Foundation of Korea (NRF) grant funded by the Korea government (MSIT) (No. RS-2023-00208874). S.-B. Yoo was supported by the National Research Foundation of Korea(NRF) grant funded by the Korea government (MSIT) (No. 2021R1F1A1062436).


\section{Preliminary}\label{collect}
In this section, we introduce a general theory about a geometric structure of a smooth projective variety with a torus action. Also we collect the well-known facts about the multiple structure on Cohen-Macaulay curves.
\subsection{Bia\l ynicki-Birula (BB) Theorem}
Let $X$ be a smooth projective variety with a $\CC^*$-action. Then the $\CC^*$-fixed locus of $X$ decomposes into
\[
X^{\CC^*}=\bigsqcup\limits_{i} Y_i
\]
such that each component $Y_i$ is connected. Note that $Y_i$ is smooth (\cite{Ive72}). For each $Y_i$, the $\CC^*$-action on the tangent bundle $TX {\mid}_{Y_i}$ provides a decomposition as
\[
TX{\mid}_{Y_i}=T^{+}\oplus T^{0}\oplus T^{-}
\]
where $T^+$, $T^{0}$ and $T^{-}$ are the subbundles of $TX{\mid}_{Y_i}$ such that the group $\CC^*$ acts with positive, zero and negative weights respectively. Under the local linearization, $T^0\cong TY_i$ and
\[
T^+\oplus T^{-}=N_{Y_i/X}=N^{+}(Y_i)\oplus N^{-}(Y_i)
\]
is the decomposition of the normal bundle $N_{Y_i/X}$ of $Y_i$ in $X$.

A fundamental result in a theory of $\CC^*$-action on $X$ has been provided by A. Bia{\l}ynicki-Birula. Let
\[
X^{+}(Y_i)=\left\{x\in X\;\mid ~\lim_{t\to 0} t\cdot x \in Y_i\right\}\quad  \text{and}\quad X^{-}(Y_i)=\left\{x\in X \;\mid ~\lim_{t\to \infty} t\cdot x \in Y_i\right\}.
\]

\begin{theorem}[\protect{\cite{Bir73}}]\label{BBthm}
Under the above assumptions and notations,
\begin{enumerate}
\item $X=\bigsqcup\limits_{i} X^{\pm}(Y_i)$;
\item for each connected component $Y_i$, there are $\CC^*$-equivariant isomorphism $X^{\pm}(Y_i)\cong N^{\pm}(Y_i)$ over $Y_i$ where $N^{\pm}(Y_i)\lr Y_i$ is a Zariski locally trivial fiberation with the affine space $\CC^{\mu^{\pm}(Y_i)}$ of dimension $\mu^{\pm}(Y_i):=\rank N^{\pm}(Y_i)$.
\end{enumerate}
\end{theorem}

Throughout this article, we denote a smooth quadric threefold simply by $Q_3=Q$. 

\begin{proposition}
The moduli space $\bS_d(Q) (\cong \bH_d(Q))$ is smooth for $d\leq 3$.
\end{proposition}
\begin{proof}
The case $d=1$ is easily verified by the normal bundle sequence of $L\subset Q\subset \PP^4$. The cases $d=2$ and $d=3$ have been proved in \cite[Theorem 3.7 and Proposition 4.13]{CHK12}.
\end{proof}


\subsection{Non-reduced cubic curves in a quadric threefold}\label{2.2}
Note that $Q$ is defined by a quadratic polynomial and does not contain any plane. Hence each curve $C$ with Hilbert polynomial $\chi(\cO_C(m))=3m+1$ is Cohen-Macaulay (CM). Let us start by recalling the list of the non-reduced CM curve $C$ in a quadric threefold $Q$; see \cite{EH82} and \cite[Lemma 2.1]{Chu22}. Let $p_a(C):=\dim\rH^1(C,\cO_C)$ be the \emph{arithmetic genus} of the curve $C$.

\begin{enumerate}[(i)]
\item (Triple thickness) The structure sheaf $\cO_C$ with $[C]=3[L]$ fits into the non-split extension
\[
\ses{\cO_L(-1)\oplus \cO_L(-1)}{\cO_C}{\cO_L}.
\]
Moreover, such $C$'s are parameterized by the GIT-quotient:
\begin{equation}\label{gitthick}
\PP(\Ext_{Q}^1(\cO_L,\cO_L(-1)\oplus \cO_L(-1)))^{\mathrm{s}}\git \mathrm{SL}(2) \cong \mathrm{Gr}(2, \mathrm{dim} \Ext_{Q}^1(\cO_L, \cO_L(-1))).
\end{equation}
\item (Triple line lying on a quadric cone) The structure sheaf $\cO_C$ with $[C]=3[L]$ fits into the non-split extension
\[
\ses{\cO_{L}(-1)}{\cO_C}{\cO_{L^2}},
\]
where $p_a(L^2)=0$.
\item (Pair of lines) The structure sheaf $\cO_C$ with $[C]=2[L]+[L']$ fits into the non-split extension
\[
\ses{\cO_{L}(-1)}{\cO_C}{\cO_{L\cup L'}},
\]
where $p_a(L^2)=0$ or $-1$.
\end{enumerate}

\begin{remark}
More generally, let $G$ be a stable sheaf on a projective space $\PP^n$ with $\chi(G(m))=(d-1)m + 1$ and $L$ be a line. Then every sheaf $F$ fitting into the non-split short exact sequence
\[\ses{\cO_L(-1)}{F}{G}\]
is stable; see \cite[Lemma 4.7]{CCM14}.
\end{remark}

\subsection{Deformation theory}
We address an exact sequence which we will use later.
\begin{lemma}\label{nestedseq}
Let $Y \stackrel{i}{\hookrightarrow} X$ be a smooth, closed subvariety of the smooth variety $X$.
If $F$ and $G \in \mathrm{Coh}(Y)$, then there is an exact sequence
\begin{equation}\label{thomas2}
\begin{split}
0 \to \Ext^1_{Y}(F, G) \to \Ext^1_{X}(i_*F, i_*G) &\to \Hom_{Y}(F,G\otimes N_{Y/X}) \\
&\to \Ext^2_{Y}(F, G) \to \Ext^2_{X}(i_*F, i_*G).
\end{split}
\end{equation}
\end{lemma}

\begin{proof}
This is the base change spectral sequence in \cite[Theorem 12.1]{MCC00}.
\end{proof}
\section{Application of the BB-theorem on a quadric threefold}\label{bb}

In this section we find the fixed rational curves and its weights in $Q$ under the $\CC^*$-action which inherits from the $\SL_2:=\text{SL}(2,\mathbb{C})$-representation. We will use the notation $\fsl_{2}$ for the associated Lie algebra to $\SL_2$. Let us regard the quadric threefold $Q$ as a hyperplane section of $\Gr(2,4)$. The $\CC^*$-actions on $Q$ and curves in there come from an action on its related vector space.

Let $V_{d}=\Sym^d\left(\CC^{2}\right)$ with the $\SL_{2}$-action induced from the left multiplication of $\SL_2$ on $\CC^{2}$.
Then the maximal torus $\{\text{diag}(t^{-1},t) \mid t\in \CC^*\}\cong \CC^{*}$ acts on the basis vectors 
\[
\left\{v_{d}=x^{d},\quad v_{d-2}=x^{d-1}y,\quad \dots \quad, v_{-1}=xy^{2-d},\quad v_{-d}=y^{d}\right\}
\]
with weights $(d, d-2, \dots, 2-d, -d)$. The infinitesimal $\fsl_{2}$-action on $V_{d}$ is given by 
\[
e=x\partial_{y},\quad f=y\partial_{x},\quad h=x\partial_{x}-y\partial_{y}
\]
for the standard basis $\langle e,f,h\rangle $ for $\fsl_{2}$. Indeed, we have the following equations for $0\leq i\leq d$:
\[
\left\{\hspace{.3cm}
\begin{split}
e{\cdot} v_{d-2i}&= iv_{d-2(i-1)},\\
f{\cdot} v_{d-2i}&=(d-i)v_{d-2(i+1)} ,\\
h{\cdot} v_{d-2i}&= (d-2i)v_{d-2i}.
\end{split}\right.
\]
Now fix $d=3$ and set $V=V_3$. Setting $v_{a,b}:=v_a^* \wedge v_b^* \in \wedge^2 V^*$ and $W:=\wedge^2 V$, we get the Pl\"{u}cker coordinates for $i : \Gr(2, V)\hookrightarrow\PP(W)$
\[
v_{3,1},~v_{3,-1},~v_{3,-3},~v_{1,-1},~v_{1,-3},~v_{-1,-3}
\]
with weights $(4,2,0,0,-2,-4)$ and the defining equation of $\Gr(2,V)$ is given by the Klein relation:
\begin{equation}
v_{3,1}v_{-1,-3}-v_{3,-1}v_{1,-3}+v_{3,-3}v_{1,-1}=0.
\end{equation}
Note that for $g \in \fsl_{2}$, $g(v_{i,j})=(g{\cdot} v_i^{*})\wedge v_j^{*} +v_i^{*}\wedge (g{\cdot} v_j^{*})$ for $i\neq j$. For example, $g=e$ acts on the basis of $W$ by
\begin{align*}
e{\cdot} v_{3,1}&=0, &e{\cdot} v_{3,-1}&=2v_{3,1}, &e{\cdot} v_{3,-3}&=3v_{3,-1}\\
e{\cdot} v_{1,-1}&=v_{3,-1}, &e{\cdot} v_{1,-3}&=v_{3,-3}+3v_{1,-1}, &e{\cdot} v_{-1,-3}&=2v_{1,-3}.
\end{align*}
By the action of $h\in \fsl_{2}$, each $v_i$ is an eigenvector of $V$ so that we get a decomposition
\[
V=V(3)\oplus V(1)\oplus V(-1)\oplus V(-3)
\]
such that $V(i)=\langle v_i \rangle$ for each $i\in \{3,1,-1,-3\}$. Similarly we get a decomposition for $W$:
\[
W=W(4)\oplus W(2)\oplus W(0)\oplus W(-2)\oplus W(-4),
\]
where $W(4)=\langle v_{3,1}^*\rangle$, $W(2)=\langle v_{3,-1}^*\rangle$, $W(0)=\langle v_{3,-3}^*, v_{1,-1}^*\rangle$ and $W(-2)=\langle v_{1,-3}^*\rangle$ and $W(-4)=\langle v_{-1,-3}^*\rangle$.
\begin{remark}\label{invhyper}
Let $H=\PP (W_1)\cong \PP^4$ be the $\SL_2$-invariant hyperplane of $\PP(W)$ for a five-dimensional subspace $W_1\subset W$. Then it must be defined by an equation $av_{3,-3}+bv_{1,-1}=0$ for some $a, b\in \CC$. In particular, we have 
\[
e{\cdot}(av_{3,-3}+bv_{1,-1})=(3a+b)v_{3,-1}=0. 
\]
Thus $H$ is defined by an linear equation $v_{3,-3}-3v_{1,-1}=0$. 
\end{remark}
Unless otherwise stated, let us define the smooth quadric threefold $Q$ by $Q=\Gr(2,V)\cap H$ for the hyperplane $H$ in Remark \ref{invhyper}.

\subsection{Fixed points and lines in $Q$}	
The fixed points $Q^{\CC^*}$ of $Q$ are the following four points:
\begin{equation}
\begin{split} 
p_{3,1}&=[1:0:0:0:0:0],\quad p_{3,-1}=[0:1:0:0:0:0],\\
p_{1,-3}&=[0:0:0:0:1:0],\quad p_{-1,-3}=[0:0:0:0:0:1].
\end{split}
\end{equation}
Note that the fixed points $H^{\CC^*}$ consist of $5$ points, including $Q^{\CC^*}$. From the defining equation of $H$, the fifth fixed point must be
\begin{equation}\label{5th}
q_0=[0:0:3:1:0:0].
\end{equation}
Therefore one can get a decomposition $W=W_1\oplus W_1^{\bot}$ so that 
\[
\PP(W_1)=\langle p_{3,1},p_{3,-1}, p_{1,-3}, p_{-1,-3}, q_0\rangle,
\]
where $\langle S \rangle$ denotes the linear span for a subset $S\subset \PP(W)$.
\begin{lemma}\label{fixedlines}
There exist exactly four lines on $Q$, fixed by the action of $\CC^{*}$: 
\begin{itemize}
\item [(a)] $\overline{p_{3,1}p_{3,-1}}=[s:t:0:0:0:0]$,
\item [(b)] $\overline{p_{3,1}p_{1,-3}}=[s:0:0:0:t:0]$,
\item [(c)] $\overline{p_{3,-1}p_{-1,-3}}=[0:s:0:0:0:t]$, and
\item [(d)] $\overline{p_{1,-3}p_{-1,-3}}=[0:0:0:0:s:t]$.
\end{itemize}
\end{lemma}

\begin{proof}
Let $L(\cong \PP^1)$ be a line in $Q$ fixed by the $\CC^*$-action. Then one can regard that the torus $\CC^*$ acts on $L$. But $\chi(L)=2$ and thus two points on $L$ are fixed. So we obtain the four lines as listed.
\end{proof}
From the inclusion $L\subset Q\subset \PP (W_1)$ and Lemma \ref{nestedseq}, we have an exact sequence:
\begin{equation}\label{eqa111}
0\to T_{[L]}\mathbf{S}_1(Q)\to  T_{[L]}\Gr(2,W_1)\to \mathrm{H}^0(\cO_{L}(2))\to 0
\end{equation}
for each line $L\subset Q$. The second term to the last in (\ref{thomas2}) is zero since $\mathbf{S}_1(Q)$ is smooth with $\dim\mathbf{S}_1(Q)=3$. 

\begin{remark}[Table \ref{hypertable1}]
Let us denote $\delta (L)$ by the number of negative weights of the tangent space of the moduli space $\mathbf{S}_1(Q)$ at the closed point $[L]$. One can compute the number $\delta(L)$ for each case in Lemma \ref{fixedlines}.

\noindent (a) The weights of the fixed line $L=\PP(\CC^2)$ are $(4, 2)$. From the isomorphism 
\[
T_{[L]}\Gr(2,W_1)\cong \Hom(\CC^2, W_1/\CC^2)=(\CC^2)^\vee \otimes (W_1/\CC^2)
\]
we get that the weight of $T_{[L]}\Gr(2,W_1)$ is given by $(-8,-6,-6,-4,-4,-2)$. Since $\mathrm{H}^0(\cO_{L}(2))\cong\Sym^2 (\CC^{2})^\vee$ has the weight $(-8,-6,-4)$, the weight of $T_{[L]} \mathbf{S}_1(Q)$ is $(-6,-4,-2)$. In particular, we get $\delta(L)=3$. 

\noindent (b) Similarly as in (a), we get that the weight of the fixed line $L=\PP(\CC^2)$ is $(4, -2)$ and so the weight of $T_{[L]} \mathbf{S}_1(Q)$ is $(-4,-2,2)$ and so $\delta(L)=2$. 

\noindent (c) and (d): Again as in the above, we get $\delta (L)=1$ (resp. $0$) for the case (c) (resp. (d)).
\end{remark}

\begin{table}[]
\centering
\begin{tabular}{|c|c|c|c|}
\hline
Type&Weights of&Weights of&Number of negative weights of \\
&$T_{[L]}\Gr(2,W_1)$ &$\mathrm{H}^0(\cO_{L}(2))$&$T_{[L]} \mathbf{S}_1(Q)$ \\ \hline
(a)&$-8,-6,-6,-4,-4,-2$&$-8,-6,-4$&$3$ \\ \hline
(b)&$-8,-4,-2,-2,2,4$&$-8,-2,4$&$2$ \\ \hline
(c)&$-4,-2,2,2,4,8$&$-4,2,8$&$1$ \\ \hline
(d)&$2,4,4,6,6,8$&$4,6,8$&$0$ \\
\hline
\end{tabular}
\captionsetup{justification=centering}
\caption{\label{hypertable1} Weights of $T_{[L]} \mathbf{S}_1(Q)$}
\end{table}

Let $P(X,t)=\sum\limits_{i=0}^{2\dim X}h^i(X, \CC)t^i$ be the Poincar\'e polynomial of a smooth projective variety $X$.
\begin{corollary}
The Poincar\'e polynomial of $\bS_1(Q)$ is
\[
P(\bS_1(Q),t)=1+t^2+t^4+t^6.
\]
\end{corollary}
\begin{proof}
The proof is straightforward by Theorem \ref{BBthm} and the result of Table \ref{hypertable1}.
\end{proof}

\subsection{Fixed conics in $Q$ and its weights}

\begin{lemma}\label{fixedconics}
There exist exactly $10$ conics on $Q$, fixed under the action of $\CC^{*}$. Furthermore, the defining equation of fixed conics are given by
\begin{itemize}
\item [(1-a)] $\langle v_{1,-3}, v_{-1,-3}, v_{3,-3}-3v_{1,-1}, v_{3,-3}v_{1,-1}\rangle$,
\item [(1-b)] $\langle v_{3,-1}, v_{-1,-3}, v_{3,-3}-3v_{1,-1}, v_{3,-3}v_{1,-1}\rangle$,
\item [(1-c)] $\langle v_{3,1}, v_{1,-3}, v_{3,-3}-3v_{1,-1}, v_{3,-3}v_{1,-1}\rangle$,
\item [(1-d)] $\langle v_{3,1}, v_{3,-1}, v_{3,-3}-3v_{1,-1}, v_{3,-3}v_{1,-1}\rangle$,
\item [(2-a)] $\langle v_{-1,-3}, v_{3,-3}, v_{1,-1}, v_{3,-1}v_{1,-3}\rangle$,
\item [(2-b)] $\langle v_{1,-3}, v_{3,-3}, v_{1,-1}, v_{3,1}v_{-1,-3}\rangle$,
\item [(2-c)] $\langle v_{3,-1}, v_{3,-3}, v_{1,-1}, v_{3,1}v_{-1,-3}\rangle$,
\item [(2-d)] $\langle v_{3,1}, v_{3,-3}, v_{1,-1}, v_{3,-1}v_{1,-3}\rangle$,
\item [(3-a)] $\langle v_{3,1}, v_{-1,-3}, v_{3,-3}-3v_{1,-1}, v_{3,-1}v_{1,-3}-v_{3,-3}v_{1,-1}\rangle$, and
\item [(3-b)] $\langle v_{3,-1}, v_{1,-3}, v_{3,-3}-3v_{1,-1}, v_{3,1}v_{-1,-3}+v_{3,-3}v_{1,-1}\rangle$.
\end{itemize}
\end{lemma}

\begin{figure}
  \centering
    \includegraphics[width=0.7\textwidth]{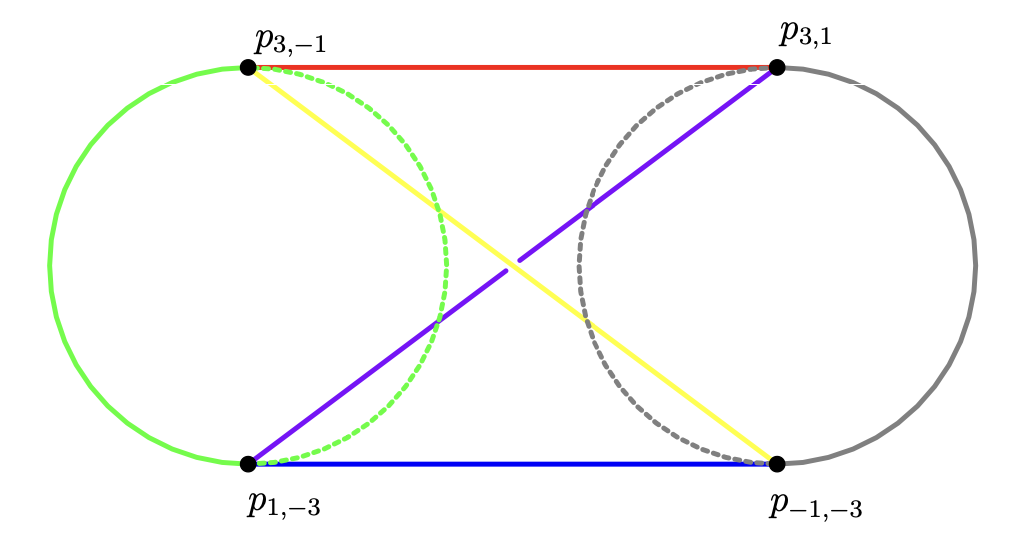}
\caption{Fixed conics in $Q$} \label{fig:M1}
\end{figure}

\begin{proof}
Let $[C]\in \mathbf{S}_2(Q)$ be a fixed conic. Then the linear span $\langle C\rangle\cong \PP^2$ is also invariant under $\CC^*$-action. The $\CC^*$-fixed planes are generated by the point $q_0$ in equation \eqref{5th} and two different points in $Q^{\CC^*}$, or three different points in $Q^{\CC^*}$. From these planes, one obtains the list.
\end{proof}
 (1-a), (1-b), (1-c) and (1-d) of the list of Lemma \ref{fixedconics} are the \emph{unique} double lines supported on (a), (b), (c) and (d) respectively in Lemma \ref{fixedlines}.  Also (2-a), (2-b), (2-c) and (2-d) in Lemma \ref{fixedconics} are  pairs of two lines supported on the fixed lines in Lemma \ref{fixedlines}. Lastly, (3-a) (resp. (3-b)) in Lemma \ref{fixedconics} are smooth conics passing through the points $p_{3,-1}$ and $p_{1,-3}$ (resp. $p_{3,1}$ and $p_{-1,-3}$). As combining these ones, we obtain its configuration (Figure \ref{fig:M1}).

Note that $\mathbf{H}_2(W_1)$ is isomorphic to $\PP(\Sym^2(\cU^\vee))$ where $\cU$ is the universal subbundle of $\Gr(3,W_1)$. Hence one can apply Lemma \ref{nestedseq} to get an exact sequence
\begin{equation}\label{conseq}
0\to T_{[C]} \mathbf{S}_2(Q)\to T_{[C]}\PP(\Sym^2(\cU^\vee))\to \mathrm{H}^0(\cO_{C}(2))\to 0
\end{equation}
for any conic $C\subset Q$. In equation \eqref{conseq}, the last term is zero since $\mathbf{S}_2(Q)$ is smooth and $\dim\mathbf{S}_2(Q)=6$.
\begin{remark}
The case (3-a) of Lemma \ref{fixedconics}: $\cU_C=\langle v_{3,-1},3v_{3,-3}+v_{1,-1},v_{1,-3}\rangle$
has weights $(2,0,-2)$ and
\begin{equation}\label{conictang}
T_{[C]}\PP(\Sym^2(\cU^\vee))=\Hom(\cU_{C},W_{1}/\cU_{C})\oplus(\Sym^2(\cU_{C}^\vee)/\CC).
\end{equation}
The first (resp. second) term in equation \eqref{conictang} has weights
$(-6,-4,-2,2,4,6)$ (resp. $(-4,-2,0,2,4)$). On the other hand, from the short exact sequence $0\to\cO_{\PP(\cU_{C})}\to\cO_{\PP(\cU_{C})}(2)\to\cO_{C}(2)\to0$, we have an exact sequence
\[
0\to \rH^{0}(\cO_{\PP(\cU_{C})})\to \rH^{0}(\cO_{\PP(\cU_{C})}(2))\to \rH^{0}(\cO_{C}(2))\to0.
\]
Since $\rH^{0}(\cO_{\PP(\cU_{C})}(2))=\Sym^{2}\rH^{0}(\cU_{C}^{*})$, $\mathrm{H}^0(\cO_{C}(2))$ has weights $(-4,-2,0,2,4)$. Therefore $T_{[C]} \mathbf{S}_2(Q)$ has weights $(-6,-4,-2,2,4,6)$. The other cases are presented in Table \ref{hypertable3}.
\end{remark}
\begin{table}[]
\begin{tabular}{|c|c|c|c|}
\hline
Type&Weights of&Weights of&Number of negative weights of \\
&$T_{[C]}\PP(\Sym^2(\cU^\vee))$ &$\mathrm{H}^0(\cO_{C}(2))$&$T_{[C]} \mathbf{S}_2(Q)$ \\ \hline
(1-a)&-8,-6,-6,-4,-4,-2,-8,-6,-4,-4,-2&-8,-6,-4,-4,-2&6 \\ \hline
(1-b)&-8,-4,-2,-2,2,4,-8,-4,-2,2,4&-8,-4,-2,2,4&4 \\ \hline
(1-c)&-4,-2,2,2,4,8,-4,-2,2,4,8&-4,-2,2,4,8&2 \\ \hline
(1-d)&2,4,4,6,6,8,2,4,4,6,8&2,4,4,6,8&0 \\ \hline
(2-a)&-8,-6,-4,-2,-2,2,-8,-6,-4,-2,4&-8,-6,-4,-2,4&5 \\ \hline
(2-b)&-6,-4,-4,-2,2,4,-8,-6,-4,2,8&-8,-6,-4,2,8&4 \\ \hline
(2-c)&-4,-2,2,4,4,6,-8,-2,4,6,8&-8,-2,4,6,8&2 \\ \hline
(2-d)&-2,2,2,4,6,8,-4,2,4,6,8&-4,2,4,6,8&1 \\ \hline
(3-a)&-6,-4,-2,2,4,6,-4,-2,0,2,4&-4,-2,0,2,4&3 \\ \hline
(3-b)&-6,-2,-2,2,2,6,-8,-4,0,4,8&-8,-4,0,4,8&3 \\
\hline
\end{tabular}
\caption{\label{hypertable3} Weights of $T_{[C]} \mathbf{S}_2(Q)$}
\end{table}
\begin{corollary}
The Poincar\'e polynomial of $\bS_2(Q)$ is
\[
P(\bS_2(Q),t)=1+t^2+2t^4+2t^6+2t^8+t^{10}+t^{12}.
\]
\end{corollary}
\begin{proof}
The proof is straightforward by Theorem \ref{BBthm} and the result of Table \ref{hypertable3}.
\end{proof}
\subsection{Fixed components of $\bS_3(Q)$: Degenerated case}
In this subsection, the fixed degenerate cubic curve is found using the fixed lines and conics (Lemma \ref{fixedlines} and Lemma \ref{fixedconics}). Furthermore, we also find fixed smooth cubics component by analyzing cubic curves lying in the quadric cone (Proposition \ref{mainprop2}). 

\begin{lemma}
Every Cohen-Macaulay curve $C\in \bS_3(Q)$ with $[C]=3[L]$ is a triple line lying on a quadric cone (Subsection \ref{2.2}).
\end{lemma}

\begin{proof}
Since $\text{dim}\Ext_{Q}^1(\cO_L,\cO_L(-1))=\mathrm{dim} \rH^0(N_{L\mid Q}(-1))=1$ for any line $L$ in $Q$ by Lemma \ref{nestedseq}, the GIT-quotient in equation \eqref{gitthick} is an empty set and so the assertion follows. 
\end{proof}

\begin{lemma}\label{fixedtriple}
Let $L$ be a $\CC^*$-invariant line and $L^2$ be the unique planar double line in $Q$ supporting $L$. Then we have
\[
\dim\Ext_Q^1(\cO_{L^2},\cO_L(-1))=2.
\]
Furthermore, the induced action on the extension space is non-trivial.
\end{lemma}

\begin{proof}
The dimension of the extension space can be easily checked by using Macaulay2 (\cite{M2}).
By taking $\Ext_Q^{\bullet}(-,\cO_L(-1))$ into the short exact sequence $\ses{\cO_L(-1)}{\cO_{L^2}}{\cO_L}$, we have a long exact sequence
\begin{equation}\label{longg}
\begin{split}
0&\cong \Hom_Q(\cO_{L^2},\cO_L(-1))\lr\Hom_Q(\cO_{L}(-1),\cO_L(-1))\stackrel{\cong}{\lr}\Ext^1_Q(\cO_{L},\cO_L(-1))\\
&\lr\Ext^1_Q(\cO_{L^2},\cO_L(-1))\lr\Ext^1_Q(\cO_{L}(-1),\cO_L(-1))\lr\cdots.
\end{split}
\end{equation}
The equality of the first term in \eqref{longg} holds due to the stability. Also the second isomorphism in \eqref{longg} holds by the dimension counting. But the last term in \eqref{longg} is isomorphic to the tangent space $T_{[L]}\bS_1(Q)$. By Table \ref{hypertable1}, there does not exist the reduplication in the space $\Ext^1_Q(\cO_{L}(-1),\cO_L(-1))$ and thus the induced action is non-trivial.
\end{proof}

\begin{lemma}\label{pairofline}\label{typeofdeg}
Let $C=L_1\cup L_2$ be a fixed conic in $Q$. Then we have
\[
\dim \Ext_{Q}^1(\cO_{C}, \cO_{L_2}(-1))=2.
\]
Furthermore, the induced action on the extension space is non-trivial.
\end{lemma}

\begin{proof}
One can prove the claim by the same method in the proof of Lemma \ref{fixedtriple}.
\end{proof}

From this discussion, one can present all of the $\CC^*$-invariant degenerated cubic curves in $\bS_3(Q)$ as follows.

\begin{proposition}\label{mainprop1}
The number of the $\CC^*$-invariant, degenerated cubic curves in $\bS_3(Q)$ is $36$.
\end{proposition}

\begin{proof}
We read the invariant cubics from Figure \ref{fig:M1}.
The reduced cubic cases are a tree of three lines or a conic attached a line. From Figure \ref{fig:M1}, one can read such cases are $4+2^3=12$. When the support of cubic is a pair of lines, then our choice is $4\times 2=8$. For each case, by Lemma \ref{typeofdeg}, there are two cases. Hence the possible case is $16$. Lastly, the support of triple line lying on a quadric cone has $4$ cases and by Lemma \ref{fixedtriple}, we have two cases for each of them. After of all, we have $12+16+8=36$.
\end{proof}

\subsection{Fixed components of $\bS_3(Q)$: Smooth case}
To find the smooth invariant cubic curves, we consider the \emph{enveloping map} 
\[
\xi:\mathbf{S}_3(Q) \rightarrow  \Gr(4,W_1)	
\]
defined by $[C]\mapsto \langle C\rangle\cong\PP^3\subset \PP(W_1)=H$. As did in the case of $\mathbf{S}_2(Q)$, we can use the five fixed points $Q^{\CC^*}\cup \{q_0\}$ to obtain the $\CC^*$-fixed points in $\mathrm{Gr}(4,W_1)$. The result is the second row of the Table \ref{hypertable}.
\begin{table}[]
\begin{tabular}{|c|cccc|}
\hline
(i)                   & \multicolumn{1}{c|}{(ii)}            & \multicolumn{1}{c|}{(iii)}          & \multicolumn{1}{c|}{(iv)}           & \multicolumn{1}{c|}{(v)}           \\ \hline
$v_{3,-3}=v_{1,-1}=0$ & \multicolumn{1}{c|}{$v_{-1,-3}=H=0$} & \multicolumn{1}{c|}{$v_{1,-3}=H=0$} & \multicolumn{1}{c|}{$v_{3,-1}=H=0$} & \multicolumn{1}{c|}{$v_{3,1}=H=0$} \\ \hline
Smooth quadric surface        &                                      & Quadric cone                        &                                     &                                    \\ \hline
\end{tabular}
\caption{\label{hypertable} Fixed hyperplane section.}
\end{table}
\noindent Let us find the fixed twisted cubics in $\xi^{-1}(H_0)\subset \mathbf{S}_3(Q)$ for a $\CC^*$-invariant element $[H_0]\in \mathrm{Gr}(4,W_1)$. That is, we seek the $\CC^*$-invariant cubics in the quadric surface $S_0:=Q\cap H_0$. Since the cases (ii) and (iii) are symmetric to (iv) and (v) in Table \ref{hypertable}, we treat the first two cases only.\\
 
$\bullet$ \textbf{Case (i) in the Table \ref{hypertable}}
\begin{proposition}
Under the notations in the above, the $\CC^*$-fixed cubics with $S_0$ smooth, are always reducible.
\end{proposition}

\begin{proof}
We fall into the case (i) of Table \ref{hypertable}, in which we can let 
\[
[v_{3,1}:v_{3,-1}:v_{1,-3}:v_{-1,-3}]\]
be the homogenous coordinates of $\PP^3$. Then the smooth quadric surface $S_0\subset \PP^3$ is defined by $v_{3,1}v_{-1,-3}-v_{3,-1}v_{1,-3}=0$ induced from the Klein relation of $\Gr(2,V)$. Set 
\[
\PP^1\times \PP^1\rightarrow  S_0\subset \PP^3,\quad [s:t]\times [v:w]\mapsto [sv:sw:tv:tw]
\]
be the $\CC^*$-equivariant map where $\CC^*$ acts on $\PP^1\times \PP^1$ diagonally with weights $(1,-1,3,-3)$. Then we get $\mathbf{S}_3(S_0)\cong  {\mid}\cO_{S_0}(1,2){\mid}\sqcup {\mid}\cO_{S_0}(2,1){\mid}$. The first component of the latter space is isomorphic to the projective space $\PP\mathrm{H}^0(\cO_{\PP^1}(1)\boxtimes\cO_{\PP^1}(2))$, and by the K{\"u}nneth formula we get
\[
\mathrm{H}^0(\cO_{\PP^1}(1)\boxtimes\cO_{\PP^1}(2))\cong\mathrm{H}^0(\cO_{\PP^1}(1))\otimes \mathrm{H}^0(\cO_{\PP^1}(2))=\CC_{s,t}\otimes \Sym^2\CC_{v,w}.
\]
Hence the weights for ${\mid}\cO_{S_0}(1,2){\mid}$ are $(\pm7,\pm5,\pm1)$ and thus the fixed points are coordinate points;  we get a similar description for the second component. 

On the other hand, let $\{sv^2, svw, sw^2, tv^2,tvw, tw^2\}$ be a basis for $\mathrm{H}^0(\cO_{S_0}(1,2))$ and its dual basis by $\{h_0,\;h_1, \;\cdots, \;h_5\}$. Then each element in the system ${\mid}\cO_{S_0}(1,2){\mid}$ is written as
\[
F:=h_0sv^{2}+h_1svw+ h_2sw^{2}+h_3tv^{2}+h_4tvw+h_5tw^{2}.
\]
Using the homogenous coordinates of $\PP^3$
\[
v_{3,1}=sv, v_{3,-1}=sw, v_{1-3}=tv, v_{-1,-3}=tw,
\]
we can write
\[
\begin{split}
sF&=h_0(sv)^{2}+h_1s^2vw+h_2(sw)^{2}+h_3(sv)(tv)+h_4stvw+h_5(sw)(tw)\\
&=h_0v_{3,1}^{2}+h_1v_{3,1}v_{3,-1}+h_2v_{3,-1}^{2}+h_3v_{3,1}v_{1-3}+h_4v_{3,-1}v_{1-3}+h_5v_{3,-1}v_{-1,-3}; \\
tF&=h_0(sv)(tv)+h_1tsvw+h_2(sw)(tw)+h_3(tv)^{2}+h_4t^2vw+h_5(tw)^{2}\\
&=h_0v_{3,1}v_{1-3}+h_1v_{3,-1}v_{1-3}+h_2v_{3,-1}v_{-1,-3}+h_3v_{1-3}^{2}+h_4v_{1-3}v_{-1,-3}+h_5v_{-1,-3}^{2}.
\end{split}
\]
Thus the defining equation of a cubic curve in $S_0$ is given by $\langle sF, tF \rangle$. In particular, one can read that $\CC^*$-fixed cubics are not irreducible.
\end{proof}
$\bullet$ \textbf{Case (ii) and (iii) in the Table \ref{hypertable}}

In these cases, the quadric surface $S_0$ is singular. By taking the $\CC^*$-orbit of the point $[1:a:3:1:3a^{-1}:0]$ with $a\neq0$, there is one-parameter family of twisted cubics \begin{equation}\label{famtw}
\cC_a:= \left\{[u^3:a u^2v:3uv^2:uv^2:3a^{-1}v^3:0]\times(a)\mid [u:v]\in\PP^1,a\in \CC^*\right\}\subset \PP^5\times \CC_{a}^*
\end{equation}
parameterized by $a\in \CC_{a}^*$. Clearly, the family in (\ref{famtw}) is invariant under the $\CC^*$-action. In fact, such a family gives a connected component in $\mathbf{S}_3(Q)^{\CC^*}$. Note that $\CC_{a}^*$ is compactified as $\PP^1\subset \bS_3(Q)$ by adding the two flat limits: A triple line lying on quadric cone and a conic attached a line.

\begin{example}
One can easily see that the family $\cC_a$ in equation \eqref{famtw} lies in the case (ii) (i.e., $S_0=Q\cap \{v_{-1,-3}=0\}$) of Table \ref{hypertable}.
We compute the number of zero-weights of $T_{[C]}\mathbf{S}_3(Q)$ at $C:=\cC\mid_{a=1}$. Let us denote the number of zero-weights of the global section $\rH^0(F)$ of a locally free sheaf $F$ by $w_0(F)$. Let $f:\PP^1\rightarrow  Q\subset \PP(W_1)\subset \PP(W)$ be a morphism defined by
\[
[u:v]\mapsto [u^3:u^2v:3uv^2:uv^2:3v^3:0].
\]
If we put the weights $(1,-1)$ of the coordinate of the domain curve and shift the weights of the codomain by $-1$, then the map $f$ is a $\CC^*$-equivariant one. We compute $w_0(f^*T_{Q})$ of the space $\mathrm{H}^0(f^*T_{Q})$ which is the first deformation space of the graph space $\text{Map}(\PP^1, Q)$. Since the map $f$ is an embedding, we have a short exact sequence
\begin{equation}\label{eqpull}
0\to f^*T_{Q}\to f^*T_{\PP^4} \to f^*N_{Q\mid\PP^4} \to 0.
\end{equation}
But we also have an exact sequence 
\[
0\to f^*\cO_{\PP^4}\to f^*\cO_{\PP^4}(1)\otimes \mathrm{H}^0(\cO_{\PP^4}(1))^*\to f^*T_{\PP^4} \to 0
\]
by pulling back the Euler sequence along the map $f$. Therefore the number of zero-weights of $f^*T_{\PP^4}$ is $w_0(f^*T_{\PP^4})=3$. Also $w_0(f^*N_{Q\mid\PP^4})=w_0(\cO_{\PP^1}(6))=1$. By plugging into equation \eqref{eqpull}, we have $w_0(f^*T_{Q})=2$. Finally our space is the quotient space by $\mathrm{SL}_2=\text{Aut}(\PP^1)$. Since $w_0(T_{\PP^1})=1$, the number of zero-weights of $T_{[C]}\mathbf{S}_3(Q)$ is one. 
\end{example}

We prove now that smooth cubic curves fixed by $\CC^*$-action are only the given one in equation \eqref{famtw}. Our strategy is to find all of the fixed curves in the resolved surface where the resolution map is $\CC^*$-equivariant. Then one can find the fixed cubics in singular quadric $S_0$. 
The following contents before the end of this section are well studied in \cite[Section 2]{Rei95}, \cite[Section 1]{Hos90} and \cite[Section 2]{Cos06}.
The resolution of the quadric cone $S_0$ at the cone point $p$ is isomorphic to the Hirzebruch surface $\FF_2:=\PP(\cO_{\PP^1}\oplus \cO_{\PP^1}(2))$. Let $C_0$ (resp. $F$) be the canonical section class (resp. fiber class) in the divisor class $\text{Div}(\FF_2)=\langle C_0, F\rangle$. Then $S_0$ is given by the image of the complete linear system ${\mid}C_0+2F{\mid} (\cong \PP^3)$. Also the surface $\FF_2$ embeds in the complete linear system ${\mid} C_0+3F{\mid}(\cong \PP^5)$. That is, we have a commutative diagram:
\begin{equation}\label{comdia1}
\xymatrix{
\FF_{2}\ar@{^{(}->}[rr]^{\mid C_0+3F\mid}\ar[rrd]^{\mid C_0+2F\mid}&&\PP^5\ar[d]^{\pi}\\
&&S_0\subset\PP^3.
}
\end{equation}
Furthermore, the right vertical map $\pi$ in \eqref{comdia1} is a linear projection. 
\begin{remark}[\protect{cf. \cite[Lemma 4.12]{CPS19}}]\label{inverseimage}
Let $[C]\in \bS_3(S_0)$. Note that $C$ is a CM-curve. Let $p$ be the cone point of $S_0$. Let $\text{mult}_p(C)=1$. Then by the projection formula, the strict transform $\tilde{C}$ of $C$ along the map $\pi_{\FF_2}:\FF_2\lr S_0$ lies in the linear system $[\tilde{C}]\in {\mid} C_0+3F{\mid}$.
In fact, $\tilde{C}\cdot (C_0+2F)=C\cdot H_{\PP^3}=3$ by the projection formula. Also, $\tilde{C}\cdot C_0=1$. Hence $\tilde{C}=C_0+3F$.
Let $\text{mult}_p(C)>1$ and $\tilde{C}=\pi_{\FF_2}^{-1}(C)$ be the scheme theoretic inverse image  of $C$ along the map $\pi_{\FF_2}$. Then $[\tilde{C}]=[C_0]+[\tilde{C}_{s}]$ where $\tilde{C}_{s}$ is the strict transform of $C$ such that $[\tilde{C}_{s}]=3[F]$. Therefore $\tilde{C}$ lies in the linear system ${\mid} C_0+3F{\mid}$.
\end{remark}
We interpret $\FF_2$ and $S_0$ as \emph{scrolls} to make all relevant maps of the diagram \eqref{comdia1} a $\CC^*$-equivariant one. By definition, the \emph{rational normal scroll} $S(p,q)$ is the \emph{join variety} of the rational normal curves of degree $p$ and $q$; here, we allow $p$ or $q$ to be zero. In our case, $\FF_2\cong S(1,3)$ and $S_0\cong S(0,2)$. Then it is well-known that the map $\pi$ in \eqref{comdia1} is a successive projection from a point on the ruling line and cubic curve.
Let us present $\pi$ in \eqref{comdia1} as homogeneous coordinates. To do this, let us define $S(1,3)$ and $S(0,2)$ as the closure of the image of $4$-dimensional complex torus $(\CC^*)^4$ as below.
\begin{equation}\label{f2f2}
\xymatrix{
(\CC^*)^4\ar[rrrrrrr]^{(t_0,t_1,u_0,u_1)\mapsto[t_0^3u_0:t_0^2t_1u_0:t_0t_1^2u_0:t_1^3u_0:t_0u_1:t_1u_1]}\ar[rrrrrrrd]^{(t_0,t_1,u_0,u_1)\mapsto[t_0^3u_0:t_0^2t_1u_0:t_0t_1^2u_0:t_0u_1]}&&&&&&&\PP^5\ar[d]\\
&&&&&&&\PP^3
}
\end{equation}

Then the scroll $S(1,3)$ (resp. $S(0,2)$) is defined by the maximal minors of the \emph{catalecticant} matrix
\begin{equation}\label{cateq}
\begin{bmatrix}
z_0&z_1&z_2&z_4\\
z_1&z_2&z_3&z_5
\end{bmatrix}\;\;\; (\text{resp.} \begin{bmatrix}
z_0&z_1\\
z_1&z_2
\end{bmatrix}).
\end{equation}
Under this setting, the map $\pi$ in \eqref {f2f2} is
\[\pi([z_0:z_1:z_2:z_3:z_4:z_5])=[z_0:z_1:z_2:z_4].\]
Now we are ready to state our main proposition.
\begin{proposition}\label{mainprop2}
In the case (ii) and (iii) (and thus (iv) and (v)) of Table \ref{hypertable}, the fixed cubic curves are classified by the following.
\begin{enumerate}
\item If $S_0=Q\cap \{v_{-1,-3}=0\}$, then the $\CC^*$-fixed, irreducible cubic curve in a quadric cone $S_0$ is defined by the family $\cC_a$ in equation \eqref{famtw}.
\item If $S_0=Q\cap \{v_{1,-3}=0\}$, then every $\CC^*$-fixed cubic curve in a quadric cone $S_0$ is isolated and degenerated one.
\end{enumerate}
\end{proposition}
\begin{proof}
Case (1): Since $S_0$ is defined by $v_{3,-1}v_{1,-3}-3v_{1,-1}^2=0$, we may assume that the weights of $z_0$, $z_1$, $z_2$ and $z_4$ are $2$, $0$, $-2$ and $4$ (after rescaling the coordinate $z_1$). If we let the weights $z_3$ and $z_5$ by $-4$ and $2$ to be invariant of equations \eqref{cateq}, then the map $\pi$ is a $\CC^*$-equivariant one. This induces a $\CC^*$-action on the complete linear system ${\mid} C_0+3F{\mid}$ which is regarded as the space of quartic curves in $\FF_2$. Hence the repeat of the weight $2$ gives us one-parameter family of quartic curves in $\FF_2$ and thus rational cubic curves in $S_0$ (Remark \ref{inverseimage}). Note that the other case is isolated and the corresponding curve is a degenerated one by substituting in equation \eqref{cateq}. Obviously, the family $\cC_a$ in \eqref{famtw} is invariant under $\CC^*$-action and thus we proved the claim.

Case (2): The similar computation as in the case (1) shows that if we let the weights of $z_0, \cdots, z_4$ and $z_5$ by $4$, $0$, $-4$, $-8$, $2$ and $-2$, then $\pi$ is $\CC^*$-equivariant. Thus $\CC^*$-fixed curves in $\FF_2$ are isolated and degenerated one. Thus the same thing holds on $S_0$.
\end{proof}
It seems possible to compute the weights by using the result in \cite{ES96} and Lemma \ref{nestedseq}. Since we explicitly describe the moduli space $\bS_3(Q)$ in terms of a projective bundle in the following section, we omit it.
\section{Cubic curves space as a projective bundle}\label{sub:defq} 
In this section, we describe the moduli space $\bS_3(Q)$ as a projective bundle over a Grassmannian (Proposition \ref{mainprop4}). As a corollary, we give the cohomology ring structure of $\bS_3(Q)$ (Corollary \ref{maincor}). 
\subsection{Sheaf theoretic description of $\bS_1(Q)$} Let us recall the global description of $\bS_1(Q)$. Let $L$ be a line in $Q$. It is known that the locally free resolution of the ideal sheaf $\cI_{L,Q}$ is given by
\begin{equation}\label{lineQ}
0\to \cO_Q(-1) \to \cU_Q \to \cI_{L, Q} \to 0,
\end{equation}
where $\cU_{Q}$ is the restriction of the universal subbundle of Grassmannian $\Gr(2,4)$; see \cite[Section 2.1.1]{Fae14}. Also, every line in $Q$ arises in this fashion, so called the Hartshorne-Serre correspondence (\cite{Arr07}). Thus we get the following. 

\begin{proposition}\label{lineinq}
Let $\mathbf{H}_1(Q)$ be the Hilbert scheme of lines in $Q$. Then, 
\[
\mathbf{H}_1(Q)\cong \PP \Hom(\cO_{Q}(-1),\cU_{Q})\cong \PP\mathrm{H}^0(\cU_{Q}(1))\cong \PP^3.
\]
\end{proposition}
\subsection{Proof of Theorem \ref{mainthm}}
The moduli space of representations of a Kronecker quiver parametrizes the isomorphism classes of stable sheaf homomorphisms 
\begin{equation}\label{res1}
\cO_Q^{\oplus 2}\longrightarrow \cU_{Q}(1)
\end{equation}
up to the natural action of the automorphism group $\CC^{*}\times \GL_{2}/\CC^{*} \cong \GL_{2}$. For two vector spaces $E$ and $F$ of dimension $2$ and $1$ respectively and $V^{*} := \mathrm{H}^{0}(Q, \cU_Q(1))$, the moduli space is constructed as $\mathbf{G} := \Hom(E, V^{*}\otimes F)\git \GL_{2} \cong V^{*}\otimes F \otimes E^{*}\git \GL_{2}$ with an appropriate linearization; see \cite{Kin94}. Note that since the $\GL_{2}$ acts as a row operation on the space of $(2 \times 4)$-matrices, and thus $\mathbf{G} \cong\mathrm{Gr}(2, 4)$. 

Let $p_{1} : \mathbf{G} \times Q \to \mathbf{G}$ and $p_{2} : \mathbf{G} \times Q \to Q$ be the natural projections, and write $\cA \boxtimes \cB:=p_1^*\cA \otimes p_2^*\cB$ for $\cA \in\text{Coh}(\mathbf{G})$ and $\cB \in\text{Coh}(Q)$. If $\cU_{\mathbf{G}}$ is the universal subbundle over $\mathbf{G}$, then there is a \emph{universal morphism}
\begin{equation}\label{univmor}
\phi \: : \: \cU_{\mathbf{G}}\boxtimes \cO_Q \longrightarrow \cO_{\mathbf{G}}\boxtimes \cU_Q(1)~;
\end{equation}
see \cite{Kin94}. Set $\cJ:=\mathrm{coker}(\phi)$ and denote $\cJ_s:=\cJ{\mid}_{\{s\}\times Q}$ for each point $s\in \mathbf{G}$. 

\begin{proposition}\label{keyprop}
For the cokernel sheaf $\cJ$ of the map in \eqref{univmor}, we have the following. 
\begin{enumerate}
\item For each point $s\in \mathbf{G}$, the restriction $\cJ_{s}$ is isomorphic to a twisted ideal sheaf $\cJ_{s}\cong \cI_{L,S}(2)$ for some hyperplane section $S:=Q\cap H$ and a line $L\subset S$.
\item $\cJ$ is flat over $\mathbf{G}$ and thus the universal morphism $\phi$ in \eqref{univmor} is injective.
\end{enumerate}
\end{proposition}

\begin{proof}
If $\cK:=\text{ker}(\phi_s)$ is non-zero, then it is a reflexive sheaf of rank one on $\{s\}\times Q\cong Q$. By the semistability of $\cO_{Q}^{\oplus 2}$ and $\cU_{Q}(1)$, the sheaf $\cK$ is isomorphic to a line bundle $\cK\cong \cO_{Q}(k)$ for some $k\in \ZZ$ and the slope condition of $\text{Im}(\phi_s)$ gives $0\leq -k \leq \frac{1}{2}$. Thus we get $k=0$ and so $\phi_s$ is not stable, i.e., $\mathrm{rank} (\mathrm{H}^0(\phi_s))=1$. Therefore $\phi_s$ is injective.

Now, for the inclusion $i_1:\cO_{Q}\rightarrow \cO_{Q}^{\oplus 2}$ into the first component, we have a commutative diagram:
\begin{equation}\label{diag1}
\xymatrix{&0&0&&\\
&\cO_{Q}\ar[r]^{\bar{\phi}_s}\ar[u]&\cI_{L,Q}(1)\ar[u]&&\\
0\ar[r]&\cO_{Q}^{\oplus 2}\ar[r]^{\phi_s}\ar[u]&\cU_{Q}(1)\ar[u]\ar[r]&\text{coker}(\phi_s)\cong \cJ_s\ar[r]&0\\
&\cO_{Q}\ar[u]_{i_1}\ar@{=}[r]&\cO_{Q}\ar[u]_{\xi}&&\\
&0\ar[u]&\ar[u]0.&&}
\end{equation}
Here, $\xi:=\phi_s\circ i_1:\cO_Q \rightarrow \cU_Q(1)$ is the composition map. Then we get $\text{coker}(\xi)\cong \cI_{L,Q}(1)$ for some line $L\subset Q$ by Proposition \ref{lineinq}. The composite map $i\circ \bar{\phi}_s : \cO_{Q}\stackrel{\bar{\phi}}{\rightarrow} \cI_{L,Q}(1) \stackrel{i}{\hookrightarrow} \cO_{Q}(1)$ is not a zero-map and thus it is injective. Hence we get $\text{coker}(i\circ \bar{\phi}_s)\cong \cO_{Q\cap H}(1)$ for some hyperplane $H\subset \PP^4$. Let $S:=Q\cap H$. Then the top row in the diagram \eqref{diag1} becomes
\[
0\to \cO_Q \cong \cI_{S,Q}(1)\stackrel{\bar{\phi}_s}{\longrightarrow} \cI_{L,Q}(1)\to \cI_{L,S}(1)\to 0
\]
and so we get $\cJ_s=\cI_{L,S}(1)$, proving the assertion (1). Now $\cJ_s$ has a constant Hilbert polynomial by the result of (1) and thus  $\cJ$ is flat over $\mathbf{G}$, which confirms (2).
\end{proof}

\begin{remark}\label{linept}
By the proof of (1) of Proposition \ref{keyprop}, each hyperplane section $Q\cap H$ is parameterized by $\mathbf{G}$. Hence one can consider the universal family of hyperplane sections in $\mathbf{G}\times Q$. Furthermore, $\mathbf{G}\cong \mathrm{Gr}(2,4)$ is the space of lines in $\PP\mathrm{H}^0(\cU_{Q}(1))=\PP^3$. On the other hand, the latter space $\PP^3$ is the Fano scheme of lines in $Q$; see Proposition \ref{lineinq}. Let $\cC \subset \PP^3 \times Q$ be the universal lines over $\PP^3$ with $p: \cC \rightarrow \PP^3$ and $q: \cC \rightarrow Q$ the projection maps. Then it can be easily checked that the transform $p(q^{-1}(S\cap H))$ of the hyperplane section $Q\cap H$ becomes a line $L\subset \PP^3$, i.e., $[L]\in \mathbf{G}$.
\end{remark}

From the proof of Proposition \ref{keyprop}, we obtain the exact sequence
\begin{equation}\label{exactt}
0\to \cU_{\mathbf{G}}\boxtimes \cO_{Q} \to \cO_{\mathbf{G}} \boxtimes\cU_{Q}(1) \to \cJ \to 0.
\end{equation}
By applying the functor $R^{\bullet} p_{1,*}((-)\boxtimes \cO_{Q}(1))$ to the exact sequence \eqref{exactt}, we have
\begin{equation}\label{univseq}
0\to \cU_{\mathbf{G}} \otimes \mathrm{H}^0 (\cO_{Q}(1))\to    \cO_{\mathbf{G}} \otimes \mathrm{H}^0 (\cU_{Q}(2)) \to p_{1,*} (\cJ\boxtimes \cO_Q(1))\to 0,
\end{equation}
since $\mathrm{H}^1(\cO_{Q}(1)) = 0$. Since $\mathrm{h}^0 (\cO_{Q}(1)) = 5$ and $\mathrm{h}^0 (\cU_{Q}(2)) = 16$, we see that the direct image sheaf $p_{1,*} (\cJ\boxtimes \cO_Q(1))$ in \eqref{univseq} is a vector bundle of rank $6$ on $\mathbf{G}$.

\begin{proposition}\label{mainprop4}
The moduli space $\mathbf{S}_3 (Q)$ is isomorphic to the projective bundle $\PP(p_{1,*} (\cJ\boxtimes \cO_Q(1)))$.
\end{proposition}

\begin{proof}
Let $\cG:=p_{1,*} (\cJ\boxtimes \cO_Q(1))$ and $\pi: \PP(\cG)\rightarrow \mathbf{G}$ be the bundle morphism. Then, there exists a commutative diagram:
\[
\xymatrix{\mathbf{G}\times Q\ar[d]^{p_1}&\PP(\cG)\times Q\ar[l]_{\pi\times i}\ar[d]^{\bar{p}_1}\\
\mathbf{G}&\PP(\cG)\ar[l]^{\pi},}
\]
where $\bar{p}_1:\PP(\cG)\times Q\rightarrow \PP(\cG)$ is the projection map. Let
\[
c~:~ \bar{p}_1^*\cO_{\PP(\cG)}(-1)\to (\pi \times i)^*(\cJ\boxtimes \cO_Q(1))
\]
be the composition of the pullback of the tautological map
 \begin{equation}
\cO_{\PP(\cG)}(-1)\rightarrow  \pi^*\cG\cong \bar{p}_{1,*}(\pi\times i)^*\left(\cJ\boxtimes \cO_Q(1)\right)
\end{equation}
and the natural map 
\[
\bar{p}_{1}^{*}\bar{p}_{1,*}(\pi\times i)^*(\cJ\boxtimes \cO_Q(1))\rightarrow  (\pi\times i)^*(\cJ\boxtimes \cO_Q(1)).
\]
Let $\bar{\cZ}:=(\pi\times i)^*\cZ$ be the pull-back of the universal hyperplane section $\cZ\subset \mathbf{G}\times Q$ parameterized by $\mathbf{G}$; see Remark \ref{linept}. We claim that the local extension space $\mathcal{E}xt_{\bar{\cZ}}^1(\text{coker}(c),\cO_{\bar{\cZ}})$ is a flat family of stable sheaves over $\PP(\cG)$ with Hilbert polynomial $3m+1$. Note that it is enough to check the claim fiberwise. Over a point $x\in \PP(\cG)$, the exact sequence 
\[
0\lr \text{Im}(c) \lr (\pi\times i)^*(\cJ \boxtimes \cO_Q(1)) \lr \text{coker}(c) \lr 0
\]
becomes a short exact sequence
\[
0\to \cO_S \to \cI_{L,S}(2)\to  \cT \to 0,
\]
where $\cT:=\text{coker}(c)_{x}$. By applying the dual functor $\mathcal{H}om_{S} (-, \cO_S)$ to the exact sequence, we have
\[
\begin{split}
0\cong \mathcal{H}om_{S} (\cT, \cO_S)\to \mathcal{H}om_{S} (\cI_{L,S}(2), \cO_S) &\to \mathcal{H}om_{S} (\cO_S, \cO_S)\cong \cO_S \\
&\to \mathcal{E}xt_{S}^1 (\cT, \cO_S)\to \mathcal{E}xt_{S}^1  (\cI_{L,S}(2), \cO_S)\cong 0. 
\end{split}
\]
The first term vanishes because $\cT$ is one-dimensional. Also, the last vanishing is a special case of a more general vanishing $\mathcal{E}xt_{S}^{i \ge 1}  (\cI_{L,S}(2), \cO_S)\cong 0$, which is obvious when $S$ is smooth. If $S$ is not smooth, one can check this by using Macaulay2 (\cite{M2}). 
By computing the Hilbert polynomial of $\mathcal{H}om_{S} (\cI_{L,S}(2), \cO_S)$ and $\cO_S$, we conclude that $\mathcal{E} xt_{S}^1 (\cT, \cO_S)=\cO_C$ for some twisted cubic curve $C\subset Q$. Hence by the universal property of moduli space $\mathbf{S}_3(Q)$, there exists a morphism 
\begin{equation}\label{univmap}
\Phi:\PP(\cG)\longrightarrow  \mathbf{S}_3(Q)
\end{equation}
induced by the flat family $\mathcal{E}xt_{\bar{\cZ}}^1(\text{coker}(c),\cO_{\bar{\cZ}})$ over $\PP(\cG)$.

Lastly, we prove that the induced map $\Phi$ in \eqref{univmap} is an isomorphism. By Lemma \ref{pic} below and Zariski main theorem, it is enough to check that the map $\Phi$ is generically one-to-one. Let us choose a smooth twisted cubic $C\subset Q$ such that, for the linear span $H_0:=\langle C\rangle \cong\PP^3$, the hyperplane section $Q\cap H_0=:S_0$ is a smooth quadric surface. In $C\subset S_0$, the curve class $[C]$ in $S_0$ is automatically determined. Hence the inverse image $\Phi^{-1}([\cO_C])$ is a unique point in $\PP(\cG)$.
\end{proof}

\begin{lemma}\label{pic}
The rank of the Picard group $\mathrm{Pic}(\bS_3(Q))$ is $2$.
\end{lemma}
\begin{proof}
Following the blowing-up/down diagram in \cite{CHK12}, we know that the rank of Picard group of $\bS_3(Q)$ is the same as that of $\bM_3(Q)$. Here $\bM_3(Q)$ is the moduli space of stable maps of degree $3$ in $Q$. But the Picard group of $\bM_3(Q)$ is generated by the boundary divisor of reducible curves and the locus of stable maps whose
images meet a fixed line in $Q$ (compare with the result in \cite{Opr05}).
\end{proof}

\begin{corollary}\label{maincor}
The cohomology ring of $\bS_3(Q)$ is given by
\[\begin{split}
&\rH^*(\bS_3(Q),\QQ)\cong\QQ[c_1,c_2, h]/I,\\
I=\langle c_1^3-2c_1c_2,c_1^4-3c_1&^2c_2+c_2^2, h^6-5c_1h^5+(15c_1^2-5c_2)h^4-40c_1c_2h^3+50c_2^2h^2\rangle
\end{split}
\]
where $c_1=c_1(\cU_{\bG})$, $c_2=c_2(\cU_{\bG})$ and $h=c_1(\cO_{\PP}(1))$ is the hyperplane class of the projective bundle in Proposition \ref{mainprop4}.
\end{corollary}

\begin{proof}
By the exact sequence \eqref{univseq} and the presentation of the cohomology ring of $\Gr(2,4)$ in \cite[Theorem 5.26]{EH16}, one can obtain the result.
\end{proof}

\bibliographystyle{alpha}
\newcommand{\etalchar}[1]{$^{#1}$}

\end{document}